\documentclass[reqno,a4paper,12pt]{amsart} 

\usepackage{amsmath,amscd,amsfonts,amssymb}
\usepackage{mathrsfs,dsfont}

\usepackage{tikz}
\usetikzlibrary{patterns}
\usepackage{float}
\usepackage{appendix}
\usepackage{caption}
\usepackage{subcaption}
\usepackage{enumerate}

\numberwithin{equation}{section}
\numberwithin{figure}{section}

\def\R{\mathbb{R}}

\def\dH{\dim_{\mathcal{H}}}

\renewcommand\leq{\leqslant}
\renewcommand\geq{\geqslant}

\theoremstyle{plain}

\newtheorem*{claim*}{Claim}
\newtheorem*{thm*}{Theorem}

\theoremstyle{definition}

\newtheorem*{definition*}{Definition}
\newtheorem*{remarks*}{Remarks}
\newtheorem*{remark*}{Remark}

\newenvironment{enumerate-math}
{\begin{enumerate}
\addtolength{\itemsep}{5pt}
}
{\end{enumerate}}

\newenvironment{enumerate-text}
{\begin{enumerate}
\addtolength{\itemsep}{5pt}
}
{\end{enumerate}}

\begin{document}

\title{The group action method and radial projection}

\author{Guo-Dong Hong} 
\address{Department of Mathematics, National Taiwan University,Taipei, 106 Taiwan}
\email{kirito21734@gmail.com}
\author{Chun-Yen Shen}
\address{Department of Mathematics, National Taiwan University, and National Center for Theoretical Sciences, Taipei, 106 Taiwan}
\email{cyshen@math.ntu.edu.tw}

%\thanks{}
\date{}


\maketitle

\begin{abstract}
    
The group action methods have been playing an important role in recent studies about the configuration 
problems inside a compact set $E$ in Euclidean spaces with given Hausdorff dimension. In this paper, we further explore the group action methods to study the radial projection problems for Salem sets.

\end{abstract}

\section{\textbf{Introduction}}

The study of various projection problems for fractal sets has been one of the most active research topics in the area of geometric measure theory and harmonic analysis. For instance, the celebrated result of Marstrand \cite{Mar54} shows that one has the following results for any Borel set $E\subset\R^2$,
\begin{itemize}
 	\item [1.] if $\dH E>1$, then $|\pi_e (E)|>0$ for almost all $e\in \mathbb{S}^1$;
 	\item [2.] if $\dH E\leq 1$, $\dH \pi_e (E)=\dH E$ for almost all $e\in \mathbb{S}^1$,
\end{itemize}
where $\pi_e (E)$ denotes the orthogonal projection of $E$ into the line with direction $e$, and  $\dH E$ denotes the Hausdorff dimension of $E$. Moreover, there also have been lots of interesting works that study the exceptional sets of the projection problem. Namely from the result above by Marstrand that we know the set of directions $e \in \mathbb{S}^1$ such that $|\pi_e (E)|=0$ when $\dH E >1$ is small in terms of Lebesgue measure. However a finer result was obtained by Falconer \cite{Fal82} that if $\dH E>1$, then \begin{equation}\label{orthogonal-exception-1}\dH\left\{e\in \mathbb{S}^1 :|\pi_e (E)|=0\right\}\leq 2-\dH E,\end{equation} and the upper bound is sharp. These results in a sense show that the studies of orthogonal projection for fractal sets have been well established in part the orthogonal projection is a linear map. However not much has been known if we consider the radial projection, a non-linear projection. Given a set $E \subset \R^d$, and $x \in \R^d$, denote the map $\pi_x : E \setminus {x} \rightarrow \mathbb{S}^{d-1}$ by
$$ \pi_x(y) = \frac{x-y}{|x-y|}.$$ \par
 Along the results of Marstrand, we also can study the following two basic questions:
\begin{itemize}
 	\item [1.] For a fixed $x$, what can we say about $H^{d-1}(\pi_x(E))$;
 	\item [2.] For a fixed $x$, what can we say about $\dH  \pi_x(E) $ in terms of $\dH E$,
\end{itemize}
where $H^{d-1}$ denotes the $d-1$ dimensional Hausdorff measure. It is clear that the Hausdorff dimension of $E$ must be taken into account in these two questions. There are already some interesting results obtained by many people, we refer the readers to the survey of Mattila \cite{PM1} for known results and related questions. Recently, T. Orponen \cite{O1} has obtained very strong results for radial projection. More precisely, it was shown in \cite{O1} that for any Borel set $E \subset \R^2$ and $E$ is not contained in a line, then there always exists a point $x \in E$ such that $\dH \pi_x(E) \geq \frac{\dH E}{2}$. In addition, a recent work by P. Shmerkin \cite{PS} proves a more general version of projection theorem and connects the projection theorem to Falconer distance problem, and the work by Bochen Liu and the second listed author \cite{LS} used sum-product estimates to study the radial projection as well as the work by O. Raz and J. Zahl \cite{RZ} makes strong connections between projection theorem and other related problems.      

The purpose of this paper is to study the following version of the radial projection that is already mentioned in the paper of Orponen \cite{O1}, called the set of unit vectors spanned by $E$, namely

\[
    \mathcal{S}(E) \equiv \{\frac{x-y}{|x-y|}\in\mathbb{S}^{d-1}:x,y \in E \, and\, x\neq y\}
\] and we aim to study when we have $H^{d-1} (\mathcal{S}(E)) >0$. It is clear that the assumption $\dH E > d-1$ is in general optimal since if the set $E$ is contained in a $d-1$ dimensional plane, then the set $\mathcal{S}(E)$ has measure zero. If fact, the result of Orponen gives a quantitative version of this result. 

\begin{itemize}
    \item \textbf{Theorem 1.1 (Orponen)}  \par
    Let $\mu \in \mathcal{M}(\mathbb{R}^{d})$ be an s-Frostman measure for some $s>d-1$ with $supp(\mu)$ not contained on a hyperplane. For $p\in(1,2)$, write
    \[
        \mathcal{S}_{p}(\mu) \equiv 
        \{x\in \mathbb{R}^{d}: \pi_{x\sharp}\mu \notin \mathcal{L}^{p}(\mathbb{S}^{d-1}) \},
    \]
    where $f_{\sharp}$ denote the push-forward under $f$. \\
    Then $\dH \, \mathcal{S}_{p}(\mu) \leq 2(d-1)-s+\delta(p)$, where $\delta(p)>0$, and $\delta(p) \rightarrow 0$ as $p \searrow 1.$
\end{itemize}

A measure $\mu \in \mathcal{M}(\mathbb{R}^{d})$ is called an s-Frostman if $\mu(B(x,r))\lesssim r^s$ for all $x\in \mathbb{R}^d$ and $r\in\mathbb{R}_{+}$, where $B(x,r)$ denote the d-dimensional ball centered at $x$ with radius $r$. It is well known as Frostman's lemma that $\dH E$ = sup$\{s:$ there exists $\mu \in \mathcal{M}(\mathbb{R}^{d})$ such that $\mu$ is an s-Frostman supported on $E$ $\}$, where $E\subset \mathbb{R}^d$ is a Borel set. One can see \cite{PM2} for more details.
\par
Moreover, the radial projection not only has close connections to Falconer distance problem \cite{PS}, but it has another connections to angle problems that are already studied in the work of Harangi, Keleti, Kiss, Maga, Máthé, Mattila, and Strenner \cite{HKKMMMS}. Given a compact set $E \subset \R^d$, its associate angle set is defined as following : 
\[
    \mathcal{A}(E) \equiv \{\theta(x^{1},x^{2},x^{3}): x^{j} \in E,  \}
\]
where $\theta(x^{1},x^{2},x^{3})$ denote the interior angle of the triangle with vertices at $x^{1}x^{2}x^{3}$, at $x_{1}$. The work of Iosevich, Mourgoglou and Palsson \cite{IMP} has proved the following.
\begin{itemize}
    \item \textbf{Theorem 1.2 (Iosevich, Mourgoglou and Palsson)} \par
    Let $E \subset \mathbb{R}^{d}$ be compact with $\dH E>\frac{d+1}{2}$, then $H^{1}(\mathcal{A}(E))>0$.
\end{itemize} \par
Later the authors in \cite{IP} improved the result in theorem 1.2 to $\dim_H E>\frac{d}{2}$ which is also clearly that the threshold is optimal when $d=2$.It is worth to mention that in particular when $d=2$, $H^{1}(\mathcal{S}(E))>0$ iff $H^{1}(\mathcal{A}(E))>0$. Therefore an immediate consequence of Theorem 1.1 gives the following proposition.
\begin{itemize}
    \item \textbf{Proposition 1.3 } \par
    Let $E \subset \mathbb{R}^{2}$ be compact with $\dim_H E>1$, then $H^{1}(\mathcal{A}(E))>0$.
\end{itemize}

There also have been lots of interesting results that have been proved regarding the configurations inscribed in a given set $E$. We refer the readers to the papers \cite{GIT}, \cite{GIT2} and the references contained therein. Moreover, many of these results have been recently improved in some subsequent works such as \cite{L1}, \cite{L2} in which one of the key new ideas is to observe the invariance of the problems under some group actions. This observation was first studied in \cite{GILP} and further explored in the work \cite{L2} which in turn gives many new results for some configurations problems which include the classical Falconer distance problem. 
Motivated by \cite{L1} and \cite{L2}, we further explore the group action methods to the problem of radial projection which allow us to find new results regarding the positivity of $H^{d-1}(\mathcal{S}(E))$. 
We say a set $E \subset \mathbb{R}^{d}$ is a Salem set if for every $s< \dim_H E$, there exists $\mu \in \mathcal{M}(\mathbb{R}^{d})$ such that $|\widehat{\mu}(x)| \leq |x|^{-\frac{s}{2}}$ for all $x\in \mathbb{R}^{d}$. Configurations problems also have been studied in some recent works. For instance, 
the authors in \cite{CLP} have established the existence
of certain geometric configurations in Salem sets. Below are our main results.

\begin{itemize}
    \item \textbf{Theorem 1.4 } \par
    Let $E \subset \mathbb{R}^{d}$ be the disjoint union of compact sets $E_{1}$ and $E_{2}$ where both are Salem sets with $\dH E_{1}>\frac{d}{2}$ and $\dH E_{2}>\frac{d}{2}$, then $H^{d-1}(\mathcal{S}(E))>0$.
\end{itemize}
\begin{itemize}
    \item  \textbf{Theorem 1.5}  \par
 Let $E \subset \mathbb{R}^{2}$ be a Salem set with $\dH E> \alpha$, then 
  \[
  \dH \{x\in \mathbb{R}^{2}:H^{1}(\pi_{x}(E))=0\}\leq 4-2\alpha.
  \]
  In particular, if $\dim_H E>\frac{4}{3}$, then there exists $x\in E$ such that $H^{1}(\pi_{x}(E))>0$.
\end{itemize} 

As mentioned above the threshold $d-1$ is in general optimal for the positivity of $H^{d-1}$ measure. Our first result drops the exponent to $\frac{d}{2}$ as long as the set is a Salem set which in some sense says the set $E$ can not be concentrated in any lower dimensional space. Note that Theorem 1.1 is stronger than Theorem 1.5 when $\alpha>1$. But our Theorem 1.5 includes the case $\alpha\leq1$. We now turn to the proofs of our main results.

\section{\textbf{Proof of theorem 1.4}}
\subsection{Preliminaries}
Consider $\Phi: E \times F \rightarrow \mathbb{S}^{d-1}$ with $\Phi(x,y)=\frac{x-y}{|x-y|}$, where $E$ and $F$ are disjoint and compactly supported (and so have positive distance). Then $\Phi$ is invariant under translation and dilation, that is, $\Phi(x,y)=\Phi(x',y')$ iff $x=tx'+z$ and $y=ty'+z$ for some $t \in \mathbb{R}_{+}$ and $z \in \mathbb{R}^{d}$.
Let $\psi \in \mathcal{C}^{\infty}_{c}(\mathbb{R}^{d})$ with $0\leq \psi \leq 1$ and $\int_{\mathbb{R}^{d}} \psi(x)\, dx =1$. Set $\psi^{\epsilon}(x)=\epsilon^{-d}\psi(\frac{x}{\epsilon})$, and $\mu^{\epsilon}_{j}=\mu_{j}*\psi^{\epsilon}\in \mathcal{C}^{\infty}_{c}(\mathbb{R}^{d})$, where $\mu_{1},\mu_{2} \in \mathcal{M}(\mathbb{R}^{d})$ are measures supported on $E$ and $F$ respectively. Finally, define $\nu^{\epsilon}$ by the push-forward measure of $\mu_{1}^{\epsilon}\times\mu_{2}^{\epsilon}$ under $\Phi$.
Now we state Liu's theorem \cite{L2}:
\begin{itemize}
  \item  \textbf{Theorem 2.1 (Liu)}  \par
    Suppose that $\lambda_{G}$ is the right Haar measure on $G$, generated by translation and dilation. With the same notations above, 
    \[
        \int_{\mathbb{S}^{d-1}} |\nu^{\epsilon}(\omega)|^{2}\, d\omega
         \approx
         \int_{G} \phi(g) \prod_{j=1}^{2}[\int_{\mathbb{R}^{d}} \mu_{j}^{\epsilon}(x)\mu_{j}^{\epsilon}(gx) \, dx] \, d\lambda_{G}(g),
    \]
  where $\phi(g) \in \mathcal{C}_{c}(G)$ and the implicit constant above is independent in $\epsilon$. Moreover, if the right hand side is bounded above uniformly in $\epsilon$, then $H^{d-1}(\Phi(E,F))>0$.
\end{itemize} 
\begin{itemize}
  \item  \textbf{Remark 2.2}  \par
    The statement in theorem 2.1 is a little bit different from the original one, while the proof is almost the same, so we just point out the difference below. First of all, to use the coarea formula in the original proof, we need to assume that $\Phi$ is Lipschitz, which corresponds to our separation condition on both $E$ and $F$. Secondly, since $E$ and $F$ are both compact in our setting, so in the right hand side above, the function $\phi(g) \in \mathcal{C}_{c}(G)$ reflects this boundedness. 
\end{itemize} \par
Our first result is to use Liu's theorem to give a uniform upper bound of the right hand side. Before doing that, we first give a proposition.
   
\begin{itemize}
  \item  \textbf{Proposition 2.3}  \par
  The above quantity is controlled by the following one :
   \[
    \frac{1}{4} \sum_{i,j=1}^{2}
   \int_{0}^{\infty} \int_{\mathbb{S}^{d-1}} 
   |\widehat{\mu_{i}^{\epsilon}}(r\omega)|^{2}
   \int_{rI} |\widehat{\mu_{j}^{\epsilon}}(t\omega)|^{2} \, \frac{dt}{t} 
   \, d\omega \, r^{d-1} \, dr,
   \]
   where $I \subset \mathbb{R}_{+}$ is a compact interval.
\end{itemize} \par
(Proof) 
We first examine the integral inside:
\[
    \int_{\mathbb{R}^{d}} \mu_{j}^{\epsilon}(x)\mu_{j}^{\epsilon}(gx) \, dx =
    \int_{\mathbb{R}^{d}} \mu_{j}^{\epsilon}(x)\mu_{j}^{\epsilon}(tx-z) \, dx
\]
By Plancheral's identity in $x$, it is equal to  
\[
    t^{-n} \int_{\mathbb{R}^{d}}
    \widehat{\mu_{j}^{\epsilon}}(\xi) 
    \overline{\widehat{\mu_{j}^{\epsilon}}}(\frac{\xi}{t})
    e^{2\pi i z \cdot \frac{\xi}{t}} \, d\xi =
    \int_{\mathbb{R}^{d}}
    \widehat{\mu_{j}^{\epsilon}}(t\xi) 
    \overline{\widehat{\mu_{j}^{\epsilon}}}(\xi)
    e^{2\pi i z \cdot \xi} \, d\xi
\]
Then since $\phi(g) \in \mathcal{C}_{c}(G)$, 
\[
   \int_{G} \phi(g) \prod_{j=1}^{2}[\int_{\mathbb{R}^{d}} \mu_{j}^{\epsilon}(x)\mu_{j}^{\epsilon}(gx) \, dx] \, d\lambda_{G}(g)
   \leq
   \int_{I} \int_{\mathbb{R}^{d}} \prod_{j=1}^{2}
   [\int_{\mathbb{R}^{d}}
    \widehat{\mu_{j}^{\epsilon}}(t\xi) 
    \overline{\widehat{\mu_{j}^{\epsilon}}}(\xi)
    e^{2\pi i z \cdot \xi} \, d\xi] \, dz \, \frac{dt}{t}
\]
By Plancheral's identity in $z$, the integral in the right hand side is equal to
\[
    \int_{I} \int_{\mathbb{R}^{d}}
    \widehat{\mu_{1}^{\epsilon}}(t\xi) 
    \overline{\widehat{\mu_{1}^{\epsilon}}}(\xi)
    \widehat{\mu_{2}^{\epsilon}}(-t\xi) 
    \overline{\widehat{\mu_{2}^{\epsilon}}}(-\xi) \, d\xi \, \frac{dt}{t}
\]
By A.M. inequality, it is less than
\[
    \frac{1}{4} \sum_{i,j=1}^{2} \int_{\mathbb{R}^{d}}
    |\widehat{\mu_{i}^{\epsilon}}(\xi)|^{2}
   [ \int_{I} |\widehat{\mu_{j}^{\epsilon}}(t\xi)|^{2}
    \, \frac{dt}{t} ]\, d\xi
\]
\[
    =\frac{1}{4}  \sum_{i,j=1}^{2}
   \int_{0}^{\infty} \int_{\mathbb{S}^{d-1}} 
   |\widehat{\mu_{i}^{\epsilon}}(r\omega)|^{2}
   \int_{I} |\widehat{\mu_{j}^{\epsilon}}(tr\omega)|^{2} \, \frac{dt}{t} 
   \, d\omega \, r^{d-1} \, dr
\]
\[
   =\frac{1}{4}  \sum_{i,j=1}^{2}
   \int_{0}^{\infty} \int_{\mathbb{S}^{d-1}} 
   |\widehat{\mu_{i}^{\epsilon}}(r\omega)|^{2}
   \int_{rI} |\widehat{\mu_{j}^{\epsilon}}(t\omega)|^{2} \, \frac{dt}{t} 
   \, d\omega \, r^{d-1} \, dr
\]

\subsection{Main proof}
Since $E_{j}$'s are Salem sets for $j=1,2$, we may choose $\mu_{j} \in \mathcal{M}(\mathbb{R}^{d})$ such that $|\widehat{\mu_{j}}(t\omega)| \leq |t|^{-\frac{s}{2}}$, where $s\in [\frac{d}{2}, min\{\dH E_{1}, \dH E_{2}\}]$.
By proposition 2.3, it suffices to show the following integral in finite:
\[
    \int_{0}^{\infty} \int_{\mathbb{S}^{d-1}} 
   |\widehat{\mu_{i}^{\epsilon}}(r\omega)|^{2}
   \int_{rI} |\widehat{\mu_{j}^{\epsilon}}(t\omega)|^{2} \, \frac{dt}{t} 
   \, d\omega \, r^{d-1} \, dr
\] 
where $\widehat{\mu_{j}^{\epsilon}}=\widehat{\mu_{j}} \widehat{\psi^{\epsilon}}$ with $|\widehat{\psi^{\epsilon}}|\leq1$.
Then we have
\[
    \int_{rI} |\widehat{\mu_{j}^{\epsilon}}(t\omega)|^{2} \, \frac{dt}{t} \lesssim_{I} r^{-s}
\]
Combine with this estimate, we have
\[
    \lesssim
    \int_{0}^{\infty} \int_{\mathbb{S}^{d-1}} 
   |\widehat{\mu_{i}^{\epsilon}}(r\omega)|^{2}
   r^{d-s-1} \, dr \, d\omega
\]
\[
     \thicksim
     \int_{\mathbb{R}^{d}} 
     |\widehat{\mu_{i}}(\xi)|^{2} |\xi|^{-s} \, d\xi
     \thicksim
     \mathcal{I}_{d-s}(\mu_{i}) \leq  \mathcal{I}_{s}(\mu_{i})
     < \infty
\]
since $s>\frac{d}{2}$, and the proof is complete.

\section{\textbf{Proof of theorem 1.5} }
\subsection{Preliminaries} 
Consider the radial projection $\pi_{x}: \mathbb{R}^{d} \rightarrow \mathbb{S}^{d-1}$ with $\pi_{x}(y)=\frac{y-x}{|y-x|}$. Also, denote $\pi_{x}^{*}(\mu)=d\tau_{x}(\omega)$, where $\pi_{x}^{*}(\mu)$ means the push-forward measure of $\mu$:
\[
\int_{\mathbb{S}^{d-1}} f(\omega) \, d\tau_{x}(\omega)=
\int_{\mathbb{R}^{d}} f(\frac{y-x}{|y-x|}) \, d\mu(y)
\]
By mollification, it is equal to 
\[
\lim_{\epsilon \to 0}  
\int_{\mathbb{R}^{d}} f(\frac{y-x}{|y-x|}) \mu^{\epsilon}(y) \, dy 
\]
\[
    =\lim_{\epsilon \to 0} c_{d} 
    \int_{0}^{\infty} \int_{\mathbb{S}^{d-1}}
    f(\omega)\mu^{\epsilon}(x-t\omega) \, d\omega t^{d-1} \, dt
\]
\[
    =\lim_{\epsilon \to 0} c_{d} 
    \int_{\mathbb{S}^{d-1}} f(\omega)
    [\int_{0}^{\infty} \mu^{\epsilon}(x-t\omega) t^{d-1} \, dt]
    \, d\omega
\]
Therefore, 
$d\tau_{x}(\omega)=\lim_{\epsilon \to 0} c_{d}  \int_{0}^{\infty} \mu^{\epsilon}(x-t\omega) t^{d-1} \, dt$
as a distribution. \par
In order to obtain the pinned version estimate, we need to show the following integral is finite:
\[
    \int_{\mathbb{R}^{d}} 
    \int_{\mathbb{S}^{d-1}} (\int_{0}^{\infty} g(x-t\omega)t^{d-1}\, dt)^{2} \, d\omega
    \, d\lambda(x)
\]
where $g\in\mathcal{S}(\mathbb{R}^{d})$ and $d\lambda(x)$ is a given measure to do the exceptional estimates. While it is not easy to estimate this integral directly, so we first give a proposition to rewrite the integral:

\begin{itemize}
  \item  \textbf{Proposition 3.1}  
  \[
  \int_{\mathbb{S}^{d-1}} (\int_{0}^{\infty} g(x-t\omega)t^{d-1}\, dt)^{2} \, d\omega =
  \int_{\mathbb{S}^{d-1}} (\int_{0}^{\infty} 
  \widehat{g_{x}(\centerdot)|\centerdot|^{\frac{d}{2}}}(t\omega)
  t^{\frac{d}{2}-1}\, dt)^{2} \, d\omega,
  \]
  for all $g \in \mathcal{S}(\mathbb{R}^{d})$, where $g_{x}(y)   \equiv g(x-y)$.
\end{itemize} \par
(Proof) 
\[
    \int_{\mathbb{S}^{d-1}} (\int_{0}^{\infty} g(x-t\omega)t^{d-1}\, dt)^{2} \, d\omega
\]
\[
    =\int_{\mathbb{S}^{d-1}} (\int_{0}^{\infty} g(x-t\omega)t^{d-1}\, dt) (\int_{0}^{\infty} g(x-s\omega)s^{d-1}\, ds) \, d\omega
\]
\[
    =\int_{\mathbb{S}^{d-1}} \int_{0}^{\infty} g(x-t\omega)
    (\int_{0}^{\infty} g(x-s\omega)s^{d}\, \frac{ds}{s})t^{d-1}\, dt \, d\omega
\]
\[
    =\int_{\mathbb{S}^{d-1}} \int_{0}^{\infty} g(x-t\omega)
    (\int_{0}^{\infty} g(x-st\omega)(st)^{d}\, \frac{ds}{s})t^{d-1}\, dt \, d\omega
\]
By changing variable: $t^{d-1}\, dt \, d\omega  \thicksim \, dy$,
\[
    \thicksim
    \int_{\mathbb{R}^{d}} g(x-y)|y|^{d}
    (\int_{0}^{\infty} g(x-sy) s^{d-1}\, ds)\, dy
\]
\[
    =\int_{0}^{\infty}
    ( \int_{\mathbb{R}^{d}} g(x-y)|y|^{\alpha}
    g(x-sy)|y|^{d-\alpha }\, dy) s^{d-1} \, ds
\]
By Plancheral's identity in $y$,
\[
    =\int_{0}^{\infty}
    (
    \int_{\mathbb{R}^{d}} \widehat{g_{x}(\centerdot)|\centerdot|^{\alpha}}(\xi)
    s^{-(2d-\alpha)}
    \widehat{g_{x}(\centerdot)|\centerdot|^{d-\alpha}}(-s^{-1}\xi)\, d\xi
    )
    s^{d-1} \, ds
\]
\[
    =\int_{0}^{\infty}
    (
    \int_{\mathbb{R}^{d}} \widehat{g_{x}(\centerdot)|\centerdot|^{\alpha}}(\xi)
    \widehat{g_{x}(\centerdot)|\centerdot|^{d-\alpha}}(-s^{-1}\xi)\, d\xi
    )
    s^{-d+\alpha} \, \frac{ds}{s}
\]
\[
    =\int_{0}^{\infty}
    (
    \int_{\mathbb{R}^{d}} \widehat{g_{x}(\centerdot)|\centerdot|^{\alpha}}(\xi)
    \widehat{g_{x}(\centerdot)|\centerdot|^{d-\alpha}}(-s\xi)\, d\xi
    )
    s^{d-\alpha} \, \frac{ds}{s}
\]
By changing variable $d\xi  \thicksim t^{d-1}\, dt \, d\omega$,
\[
   \thicksim
    \int_{0}^{\infty} \int_{0}^{\infty} \int_{\mathbb{S}^{d-1}}
    \widehat{g_{x}(\centerdot)|\centerdot|^{\alpha}}(t\omega)
    \widehat{g_{x}(\centerdot)|\centerdot|^{d-\alpha}}(-s t\omega) \, d\omega \, t^{d-1} \, dt \, s^{d-\alpha} \, \frac{ds}{s}
\]
\[
    =\int_{\mathbb{S}^{d-1}}
    (
    \int_{0}^{\infty} \int_{0}^{\infty}
    \widehat{g_{x}(\centerdot)|\centerdot|^{\alpha}}(t\omega)
    \widehat{g_{x}(\centerdot)|\centerdot|^{d-\alpha}}(-s t\omega) s^{d-\alpha} \, \frac{ds}{s}\, t^{d-1} \, dt
    ) 
    \, d\omega
\]
\[
    =\int_{\mathbb{S}^{d-1}}
    (
    \int_{0}^{\infty} \int_{0}^{\infty}
    \widehat{g_{x}(\centerdot)|\centerdot|^{\alpha}}(t\omega)
    \widehat{g_{x}(\centerdot)|\centerdot|^{d-\alpha}}(-s \omega) (\frac{s}{t})^{d-\alpha} \, \frac{ds}{s}\, t^{d-1} \, dt
    ) 
    \, d\omega
\]
\[
    =\int_{\mathbb{S}^{d-1}}
    (
    \int_{0}^{\infty}
    \widehat{g_{x}(\centerdot)|\centerdot|^{\alpha}}(t\omega)
    t^{\alpha-1} \, dt
    )
    (
    \int_{0}^{\infty}
    \widehat{g_{x}(\centerdot)|\centerdot|^{d-\alpha}}(-s \omega) s^{d-\alpha-1} \, ds
    )
    \, d\omega
\]
\[
    =\int_{\mathbb{S}^{d-1}}
    (
    \int_{0}^{\infty}
    \widehat{g_{x}(\centerdot)|\centerdot|^{\alpha}}(t\omega)
    t^{\alpha-1} \, dt
    )
    (
    \int_{0}^{\infty}
    \widehat{g_{x}(\centerdot)|\centerdot|^{d-\alpha}}(s \omega) s^{d-\alpha-1} \, ds
    )
    \, d\omega.
\]
Finally, take $\alpha=\frac{d}{2}$, it is then equal to
\[
    \int_{\mathbb{S}^{d-1}} (\int_{0}^{\infty} 
  \widehat{g_{x}(\centerdot)|\centerdot|^{\frac{d}{2}}}(t\omega)
  t^{\frac{d}{2}-1}\, dt)^{2} \, d\omega.
\]
\begin{itemize}
  \item  \textbf{Remark 3.2}  \par
 
  If $d$ is an even number, we have $|y|^{d} \thicksim
  (\sum_{j=1}^{d} y_{j})^{\frac{d}{2}}(\sum_{j=1}^{d} y_{j})^{\frac{d}{2}}$ . Then
  \[
  \int_{\mathbb{S}^{d-1}} (\int_{0}^{\infty} g(x-t\omega)t^{d-1}\, dt)^{2} \, d\omega 
   \thicksim
  \int_{\mathbb{S}^{d-1}} (\int_{0}^{\infty} 
  \widehat{g_{x}(\centerdot)(\centerdot)^{\frac{d}{2}}}(t\omega)
  t^{\frac{d}{2}-1}\, dt)^{2} \, d\omega,
  \]
\end{itemize} \par
\begin{itemize}
  \item  \textbf{Remark 3.3}  \par
 If we modify the proof in proposition 3.1 and assume that $g$ is nonnegative, then we have
 \[
  \int_{\mathbb{S}^{d-1}} (\int_{0}^{\infty} g(x-t\omega)t^{d-1}\, dt)^{2} \, d\omega
  \]
  \[
    \lesssim
   \int_{\mathbb{S}^{d-1}}
    (
    \int_{-\infty}^{\infty}
    \widehat{g_{x}(\centerdot)|\centerdot|^{\alpha}}(t\omega)
    |t|^{\alpha-1} \, dt
    )
    (
    \int_{-\infty}^{\infty}
    \widehat{g_{x}(\centerdot)|\centerdot|^{d-\alpha}}(s \omega) |s|^{d-\alpha-1} \, ds
    )
    \, d\omega.
  \]
  In this case, the integrals inside integrate from $-\infty$ to $\infty$, which is helpful if we want to regard this as convolution later.\par
  In particular, 
  if $d$ is an even number, then
  \[
  \int_{\mathbb{S}^{d-1}} (\int_{0}^{\infty} g(x-t\omega)t^{d-1}\, dt)^{2} \, d\omega 
   \lesssim
  \int_{\mathbb{S}^{d-1}} (\int_{-\infty}^{\infty} 
  \widehat{g_{x}(\centerdot)(\centerdot)^{\frac{d}{2}}}(t\omega)
  t^{\frac{d}{2}-1}\, dt)^{2} \, d\omega,
  \]
  
\end{itemize} \par
From now on, we will focus on the case when $d=2$ and $\alpha=\frac{d}{2}$, then one can see the monomial term inside the integral vanishes. \\
We first examine the term 
\[
\widehat{g_{x}(\centerdot)(\centerdot)}(t\omega)
 \thicksim
\int_{\mathbb{R}^{2}} g_{x}(y)(y_{1}+y_{2})e^{-2\pi iy\cdot t\omega} \, dy
\]
\[
    \lesssim \sum_{j=1}^{2} \partial_{\xi_{j}} 
    \widehat{g_{x}(\centerdot)}(t\omega)=
    \sum_{j=1}^{2} \partial_{\xi_{j}}
    (
    \widehat{g}(\xi)e^{-2\pi i x \cdot \xi}
    )
    \vert_{\xi=t\omega}
\]
\[
    \thicksim \sum_{j=1}^{2}
    (
    \partial_{\xi_{j}} \widehat{g}(t\omega)
    e^{-2\pi i x \cdot t\omega} +
    x_{j} \widehat{g}(t\omega)
    e^{-2\pi i x \cdot t\omega}
    ).
\] \\
Then the integral under estimate becomes
\[
    \int_{\mathbb{R}^{2}} 
    \int_{\mathbb{S}^{1}} 
    (\int_{0}^{\infty} g(x-t\omega)\, dt)^{2} \, d\omega
    \, d\lambda(x)
\]
\[
    \lesssim \int_{\mathbb{R}^{2}}
    \int_{\mathbb{S}^{1}} (\int_{-\infty}^{\infty} 
  \widehat{g_{x}(\centerdot)(\centerdot)}(t\omega)
   \, dt)^{2} \, d\omega \, d\lambda(x)
\]
\[
    \lesssim \sum_{j=1}^{2}
    \int_{\mathbb{R}^{2}} \int_{\mathbb{S}^{1}}
    [
    (
    \int_{-\infty}^{\infty}
    \partial_{\xi_{j}} \widehat{g}(t\omega)
    e^{-2\pi i x \cdot t\omega} \, dt
    )^{2}+
    (
    \int_{-\infty}^{\infty}
    x_{j} \widehat{g}(t\omega)
    e^{-2\pi i x \cdot t\omega} \, dt
    )^{2}
    ]
    \, d\omega 
     d\lambda(x)
\]
\[
    = \sum_{j=1}^{2}
    \int_{\mathbb{S}^{1}} \int_{\mathbb{R}^{2}} 
    [
    (
    \int_{-\infty}^{\infty}
    \partial_{\xi_{j}} \widehat{g}(t\omega)
    e^{-2\pi i x \cdot t\omega} \, dt
    )^{2}+
    (
    |x_{j}|^{2}
    \int_{-\infty}^{\infty}
     \widehat{g}(t\omega)
    e^{-2\pi i x \cdot t\omega} \, dt
    )^{2}
    ]
    \,d\lambda(x)
    \, d\omega 
\]
\[
    \lesssim_{\lambda} \sum_{j=1}^{2}
    \int_{\mathbb{S}^{1}} \int_{\mathbb{R}^{2}} 
    [
    (
    \int_{-\infty}^{\infty}
    \partial_{\xi_{j}} \widehat{g}(t\omega)
    e^{-2\pi i x \cdot t\omega} \, dt
    )^{2}+
    (
    \int_{-\infty}^{\infty}
     \widehat{g}(t\omega)
    e^{-2\pi i x \cdot t\omega} \, dt
    )^{2}
    ]
    \,d\lambda(x)
    \, d\omega 
\]
if we assume that $d\lambda$ has compact support. \par
So we will at first deal with 
\[
    \int_{\mathbb{S}^{1}} \int_{\mathbb{R}^{2}} 
    (
    \int_{-\infty}^{\infty}
     \widehat{g}(t\omega)
    e^{-2\pi i x \cdot t\omega} \, dt
    )^{2}
    \,d\lambda(x)
    \, d\omega 
\]
and then show that the other term can be handled in a similar way.
\\
By mollification, we consider the following integral:
\[
     \int_{\mathbb{R}^{2}} 
    (
    \int_{-\infty}^{\infty}
     \widehat{g}(t\omega)
    e^{-2\pi i x \cdot t\omega} \, dt
    )^{2}
    h_{\epsilon}(x)
    \,  \lambda(x)
\]
where $h_{\epsilon}(x) \in \mathcal{S}(\mathbb{R}^{d})$,
\[
     =\int_{\mathbb{R}} \int_{\mathbb{R}} 
    (
    \int_{-\infty}^{\infty}
     \widehat{g}(t\omega)
    e^{-2\pi i (a\omega+b\omega^{\perp}) \cdot t\omega} \, dt
    )^{2}
    h_{\epsilon}(a\omega+b\omega^{\perp})
    \,  da \, db
\]
\[
     =\int_{\mathbb{R}} \int_{\mathbb{R}} 
    (
    \int_{-\infty}^{\infty}
     \widehat{g}(t\omega)
    e^{-2\pi i at} \, dt
    )^{2}
    h_{\epsilon}(a\omega+b\omega^{\perp})
    \,  da \, db
\]
\[
     = \int_{\mathbb{R}} [
    (
    \int_{-\infty}^{\infty}
     \widehat{g}(t\omega)
    e^{-2\pi i at} \, dt
    )^{2}
    \int_{\mathbb{R}}
    (
    h_{\epsilon}(a\omega+b\omega^{\perp})
    \, db
    )
    ]
    \,  da 
\]
By Plancheral's identity in $a$,
\[
    = \int_{\mathbb{R}} 
    (\widehat{g}*\widehat{g})(t\omega) 
    \widehat{h_{\epsilon}}(t\omega)\, dt
\]
where the convolution takes place in $\mathbb{R}$.

\subsection{Main proof}
Since $E$ is a Salem set, we can find $\mu\in\mathcal{M}(\mathbb{R}^{2})$ such that $|\widehat{\mu}(\xi)| \leq |\xi|^{-\frac{\alpha}{2}}$, where $\alpha\in [1, \dH E, ]$. 
Also, let 
$\lambda\in\mathcal{M}(\mathbb{R}^{2})$ be an s-Frostman. After mollification, by proposition 3.1 and the discussion above, it suffices to show that the following integral in finite:
\[
     \sum_{j=1}^{2}
    \int_{\mathbb{S}^{1}} \int_{\mathbb{R}^{2}} 
    [
    (
    \int_{-\infty}^{\infty}
    \partial_{\xi_{j}} \widehat{\mu^{\epsilon}}(t\omega)
    e^{-2\pi i x \cdot t\omega} \, dt
    )^{2}+
    (
    \int_{-\infty}^{\infty}
     \widehat{\mu^{\epsilon}}(t\omega)
    e^{-2\pi i x \cdot t\omega} \, dt
    )^{2}
    ]
    \,d\lambda(x)
    \, d\omega 
\]
And we first deal with the second term which reduced to the following integral:
\[
    \int_{\mathbb{S}^{1}}
     \int_{\mathbb{R}} 
    (\widehat{\mu^{\epsilon}}*\widehat{\mu^{\epsilon}})(t\omega) 
    \widehat{h_{\epsilon}}(t\omega)\, dt \, d\omega
\]
By Holder's inequality,
\[
    \leq
    (
    \int_{\mathbb{S}^{1}}   \int_{\mathbb{R}} 
    (\widehat{\mu^{\epsilon}}*\widehat{\mu^{\epsilon}})(t\omega)^{2} |t|^{-s+1} \, dt \, d\omega
    )^{\frac{1}{2}}
    (
    \int_{\mathbb{S}^{1}}   \int_{\mathbb{R}} 
    |\widehat{h_{\epsilon}}(t\omega)|^{2} |t|^{s-1} \, dt \, d\omega
    )^{\frac{1}{2}}
\]
while the second term
\[
    (
    \int_{\mathbb{S}^{1}}   \int_{\mathbb{R}} 
    |\widehat{h_{\epsilon}}(t\omega)|^{2} |t|^{s-1} \, dt \, d\omega
    )^{\frac{1}{2}}
     \thicksim
    \mathcal{I}_{s}(h_{\epsilon})^{\frac{1}{2}}
    \leq  \mathcal{I}_{s}(\lambda)^{\frac{1}{2}}
    < \infty
\]
So it remains to estimate the first term:
\[
    \int_{\mathbb{S}^{1}}   \int_{\mathbb{R}} 
    (\widehat{\mu^{\epsilon}}*\widehat{\mu^{\epsilon}})(t\omega)^{2} |t|^{-s+1} \, dt \, d\omega
\]
Note that we will show the estimate uniform in $\omega \in \mathbb{S}^{1}$, so we focus on 
\[
    \int_{\mathbb{R}} 
    (\widehat{\mu^{\epsilon}}*\widehat{\mu^{\epsilon}})(t\omega)^{2} |t|^{-s+1} \, dt
\]
We divide the integral above into two parts:
\[
    \int_{[-1,1]} 
    (\widehat{\mu^{\epsilon}}*\widehat{\mu^{\epsilon}})(t\omega)^{2} |t|^{-s+1} \, dt
    +
    \int_{\mathbb{R}\cap[-1,1]^{c}} 
    (\widehat{\mu^{\epsilon}}*\widehat{\mu^{\epsilon}})(t\omega)^{2} |t|^{-s+1} \, dt
\]
For the first integral, by Holder's inequality,
\[
    \int_{[-1,1]} 
    (\widehat{\mu^{\epsilon}}*\widehat{\mu^{\epsilon}})(t\omega)^{2} |t|^{-s+1} \, dt
    \leq
    \||\widehat{\mu^{\epsilon}}*\widehat{\mu^{\epsilon}}|^{2}\|_{p}\, \||t|^{1-s}\chi_{[-1,1]}\|_{p'}
\]
When $(1-s)p'>-1 \Leftrightarrow 2-s>\frac{1}{p}$, say $ 2-s-\delta =\frac{1}{p}$, the latter integral is finite. Then by Young's inequality, the former one becomes
\[
    \||\widehat{\mu^{\epsilon}}*\widehat{\mu^{\epsilon}}|^{2}\|_{p} \leq
    \|\widehat{\mu^{\epsilon}}\|^{4}_{\frac{4p}{2p+1}}
\]
Since $\widehat{\mu^{\epsilon}} \in \mathcal{S}(\mathbb{R}^{2})$, we only consider its decay rate at infinity to ensure the integral is finite. Also, we know that $|\widehat{\mu}^{\epsilon}(\xi)| \leq |\xi|^{-\frac{\alpha}{2}}$, so it reduces to the case when $(-\frac{\alpha}{2})\frac{4p}{2p+1}<-1
\Leftrightarrow 2\alpha+s > 4-\delta$.   \\ 
For the same integral, we can estimate in a similar way. In this case, we require $(1-s)p'<-1 \Leftrightarrow 2-s<\frac{1}{p}$, say $ 2-s+\delta =\frac{1}{p}$. In the end, it becomes $ 2\alpha+s > 4+\delta$. As a result, as long as $ 2\alpha+s > 4$, the integral is finite. \par
Finally, let's back to the other term
\[
     \int_{\mathbb{S}^{1}} \int_{\mathbb{R}^{2}} 
    (
    \int_{-\infty}^{\infty}
    \partial_{\xi_{j}} \widehat{\mu^{\epsilon}}(t\omega)
    e^{-2\pi i x \cdot t\omega} \, dt
    )^{2}
    h_{\epsilon}(x) \, d\omega=
    \int_{\mathbb{S}^{1}}
     \int_{\mathbb{R}} 
    ( \partial_{\xi_{j}}\widehat{\mu^{\epsilon}}* \partial_{\xi_{j}}\widehat{\mu^{\epsilon}})(t\omega) 
    \widehat{h_{\epsilon}}(t\omega)\, dt \, d\omega
\]
By property of convolution, we can pass differentiation from $\widehat{\mu^{\epsilon}}$ into $h_{\epsilon}$:
\[
    \int_{\mathbb{S}^{1}}
     \int_{\mathbb{R}} 
    ( \widehat{\mu^{\epsilon}}* \widehat{\mu^{\epsilon}})(t\omega) 
    \partial_{\xi_{j}}^{2}\widehat{h_{\epsilon}}(t\omega)\, dt \, d\omega
\]
While we still have the following control:
\[
    (
    \int_{\mathbb{S}^{1}}   \int_{\mathbb{R}} 
    |\partial_{\xi_{j}}^{2}\widehat{h_{\epsilon}}(t\omega)|^{2} |t|^{s-1} \, dt \, d\omega
    )^{\frac{1}{2}}
     \lesssim
    \mathcal{I}_{s}(|\cdot|^{2}h_{\epsilon}(\cdot))^{\frac{1}{2}}
    \lesssim_{\lambda}  \mathcal{I}_{s}(\lambda)^{\frac{1}{2}}
    < \infty
\]
since we assume that $\lambda$ has compact support. As a result, all the proof works in this case, and the proof is done, that is, 
\[
  \dH\{x\in \mathbb{R}^{2}:H^{1}(\pi_{x}(E))=0\}\leq 4-2\alpha.
  \]

\section{\textbf{Proof of proposition 1.3}}
\subsection{Preliminary} 
\begin{itemize}
    \item \textbf{lemma 4.1}  \par
    Suppose a compact set $E \subset \R^d$ has positive Hausdorff dimension. Let $\mu_E$ be a Frostman measure on $E$. Then we can find two subsets in $E$ such that each subset has positive $\mu_E$ measure and the distance between these two subsets is positive.
\end{itemize} \par
(Proof) 
We will use induction on the dimension $d$ so that we first claim it holds in $\R^1$. First since we are assuming $\dH E >0$, and thus without loss of generality we may assume $E$ does not have any point mass. After rescaling and translation, we may also assume $E \subset [0, 1]$. Now divide $[0, 1]$ into $[0, \frac{1}{2}]=I_1$, $[\frac{1}{2}, 1]= I_2$ and consider $E_1 = E \cap I_1$, $E_2 = E \cap I_2$. If both $E_1$ and $E_2$ have positive $\mu_E$ measure, then consider a small neighborhood $B_{\frac{1}{2}}(\delta)$
of the point $\frac{1}{2}$. Since we assume $E$ does not have any point mass, therefore we
can choose $\delta$ small enough so that the $\mu_E$ measure of $B_{\frac{1}{2}}(\delta)$ is smaller enough which gives the desired result. On the other hand,
if one of $E_1$ or $E_2$ has $\mu_E$ measure zero, say $\mu_E(E_2) = 0$. We repeat the arguments for $E_1$, and it is clear that we will end with a finite step $N$ so that $E_1 = E \cap I_N$ and $E_2 = E \cap I_{N+1}$ both have positive $\mu_E$ measure. Otherwise if at each step, there is only one subset $E_s = E \cap I_s$ that has positive $\mu_E$ measure, then $\mu_E(E) = \lim_{s \to \infty} \mu_E(E_s)$ which contradicts to our assumption that $E$ does not have point mass. Now assuming the result holds in $\R^{d-1}$ and we claim the result also holds in $\R^d$. Given a compact set $E \subset [0, 1]^d$ with $\dH(E) > 0$ and $\mu_E$ is a Frostman measure on $E$. Divide $[0, 1]^d$ into $2^d$ dyadic subcubes ${I_j}$. Let $E_j = E \cap I_j.$ By the same reasoning as above, we can assume two of these
subsets $E_j$ have positive $\mu_E$ measure. If these two subsets $E_1$ and $E_2$ have positive distance, then we are done. Otherwise, these two subsets $E_1$ and $E_2$ touch on a $d-1$ dimensional face, denoted
by $F$ or touch at a $d-2$ dimensional edge, denoted by $ L$. Now if $\mu_E(F) = 0$, then we can consider
a small $\delta$ neighborhood of $F$. By choosing $\delta$ small enough, we can remove this small neighborhood of $F$ from $E_1$ and $E_2$ to get what we want. If $\mu_E(F) > 0$, then by the induction hypothesis, we are done. 

\subsection{Main proof} 
Now given $E \subset \mathbb{R}^{2}$ with $\dH E>1$, by lemma 4.1, we can find $E_{1}\subset E$ and $E_{2}\subset E$ with positive distance such that $\dH E_{1}>1$ and $\dH E_{2}>1$. Then by theorem 1.1, $\dH \, \mathcal{S}_{p}(\mu_{E_{1}}) \leq 2-s+\delta(p) < 1+\delta(p)$, where $\mu_{E_{1}}$ is the s-Frostman on $E_{1}$ with $s\in (1, min\{\dH E_{1}, \dH E_{2}\})$. By choosing $p$ close to 1, we can let $\dH \, \mathcal{S}_{p}(\mu_{E_{1}})<1<\dH E_{2}$. Therefore, there exists $x \in E_{2}$ such that $\pi_{x\sharp}\mu_{E_{1}} \in \mathcal{L}^{p}(\mathbb{S}^{1})$, which means that $0<H^{1}(\pi_{x\sharp}(E_{1}))<H^{1}(\mathcal{S}(E))$.

\end{document}